\documentclass[11pt,reqno]{amsart}
\usepackage[T1]{fontenc}
\usepackage{lmodern}
\usepackage{amsmath,amssymb,amsthm,mathtools}
\usepackage[hmargin=1in,vmargin=0.95in]{geometry}
\usepackage{microtype}
\usepackage{enumitem}
\usepackage[colorlinks=true,linkcolor=blue!50!black,citecolor=blue!50!black,urlcolor=blue!50!black]{hyperref}
\usepackage{xcolor}

\newtheorem{theorem}{Theorem}[section]
\newtheorem{corollary}[theorem]{Corollary}
\newtheorem{proposition}[theorem]{Proposition}
\newtheorem{lemma}[theorem]{Lemma}
\theoremstyle{definition}
\newtheorem{conjecture}[theorem]{Conjecture}
\theoremstyle{remark}
\newtheorem{remark}[theorem]{Remark}
\numberwithin{equation}{section}
\DeclareMathOperator{\Rea}{Re}
\DeclareMathOperator{\Ima}{Im}

\newcommand{\PP}{P^{+}}
\newcommand{\e}{\mathrm e}
\newcommand{\dd}{\,\mathrm d}
\newcommand{\ii}{\mathrm i}
\newcommand{\chiz}{\chi_{0}}

\title[Smooth approximations to character sums and zeta sums]
{Smooth approximations to character sums and zeta sums}
\author[Z. Dong]{Zikang Dong}
\author[S. Zhao]{Shengbo Zhao}
\author[G. Zhou]{Guang-Liang Zhou}
\address{(Zikang Dong) School of mathematical sciences, Soochow University, Suzhou, China}
\address{(Shengbo Zhao) School of mathematical sciences, Tongji University, Shanghai, China}
\address{(Guang-Liang Zhou) College of Science, Nanjing Forestry University, Nanjing, China}
\hypersetup{pdfauthor={Zikang Dong, Shengbo Zhao, Guangliang Zhou}}
\date{}

\begin{document}
\begin{abstract}
We study smooth approximations to Dirichlet character sums and zeta sums.
The uniform upper-bound consequence of the Granville--Soundararajan conjecture,
when imposed on all nonprincipal characters including imprimitive ones,
implies the Riemann hypothesis for every primitive nonprincipal Dirichlet
\(L\)-function. The proof uses characters induced to moduli divisible by all
small primes. In the other direction, assuming the Riemann hypothesis for
\(L(s,\chi)\), we establish the conjectured approximation uniformly for
\(q^a\le x\le q\), where \(0<a\le1\) is fixed, with
\(y=(\log q+(\log x)^2)(\log\log q)^C\) for every fixed \(C>6\).
In this range both the character sum and its smooth restriction,
and hence their difference, are exponentially small relative to the
number of \(y\)-smooth integers coprime to \(q\).
We establish the corresponding approximation for
\(\sum_{n\le x}n^{-it}\) under the Riemann hypothesis for the zeta function,
uniformly for \(T^a\le x\le T\) and \(2T\le t\le5T\).
We also show that the full zeta-sum smooth-approximation conjecture
implies the Lindel\"of hypothesis.
\end{abstract}
\maketitle
\enlargethispage{6pt}

\section{Introduction and main results}

Let \(\chi\) be a Dirichlet character modulo \(q\), and let \(\chiz\) denote
the principal character modulo the same modulus. We write
\[
 S(x,\chi)=\sum_{n\le x}\chi(n),\qquad
 \Psi(x,y;f)=\sum_{\substack{n\le x\\\PP(n)\le y}}f(n),
\]
where \(\PP(n)\) is the largest prime factor of \(n\), with
\(\PP(1)=1\). As usual, \(\Psi(x,y)=\Psi(x,y;1)\). In particular,
\[
 \Psi(x,y;\chiz)=
 \#\{n\le x:\PP(n)\le y,\ (n,q)=1\}.
\]
All logarithms are natural. We write \(\log_2 q=\log\log q\) and
\(\log_3 q=\log\log\log q\), and consider \(q\) sufficiently large whenever
iterated logarithms occur. The letter \(p\) always denotes a prime.
We write \(\pi(u)\) for the number of primes not exceeding \(u\),
\(\omega(q)\) for the number of distinct prime divisors of \(q\),
and \(\Lambda(n)\) for the von Mangoldt function.
We use RH for the Riemann hypothesis for \(\zeta(s)\) and GRH for
the generalized Riemann hypothesis for Dirichlet \(L\)-functions.

Estimating \(S(x,\chi)\) is a basic problem in the study of
Dirichlet characters. Such estimates describe the distribution of
character values and enter questions about power residues, least
nonresidues, and values of Dirichlet \(L\)-functions. For primitive
characters, the P\'olya--Vinogradov inequality gives an upper bound
of order \(\sqrt q\log q\), and Montgomery and Vaughan
\cite{MV} obtained the conditional bound
\(\sqrt q\log_2q\) under GRH. For shorter sums, Burgess's estimates \cite{Burgess}
give cancellation for prime moduli whenever
\(x\ge q^{\frac{1}{4}+\varepsilon}\), where \(\varepsilon>0\) is fixed.
The interaction between character sums and zeros
of their \(L\)-functions has also been studied quantitatively by
Granville and Soundararajan \cite{GSzeros}.

A smooth approximation gives structural information in addition to
an upper bound for a character sum. The retained sum
\(\Psi(x,y;\chi)\) depends only on the values \(\chi(p)\) at
primes \(p\le y\). The restricted coefficients
\(\chi(n)\mathbf{1}_{\PP(n)\le y}\) have a Dirichlet generating
series given by a finite Euler product. Thus an approximation on this scale attributes a large
character sum to the behaviour of the character at small primes.
The corresponding unweighted count is governed by the Dickman
function: for each fixed \(u\ge1\),
\[
 \frac{\Psi(x,x^{1/u})}{x}\longrightarrow\rho(u),
 \qquad x\to\infty,
\]
where \(\rho(u)=1\) for \(0\le u\le1\) and
\(u\rho'(u)+\rho(u-1)=0\) for \(u>1\).
For the distribution of smooth integers and the estimates available
when \(u\) varies, we refer to Hildebrand and Tenenbaum \cite{HT}.

For character sums, the natural normalization also includes
coprimality to the modulus. The count \(\Psi(x,y;\chiz)\)
records both restrictions, and its dependence on \(q\) becomes
important when the modulus contains many small prime factors.
Granville and Soundararajan \cite[Conjecture~1]{GS} proposed the
following uniform approximation.

\begin{conjecture}[Granville--Soundararajan]\label{conj:GS}
There is a constant \(A>0\) such that, uniformly for every nonprincipal
character \(\chi\bmod q\) and every \(1\le x\le q\),
\begin{equation}\label{eq:GSconj}
 S(x,\chi)=\Psi(x,Y_A(x,q);\chi)
       +o\bigl(\Psi(x,Y_A(x,q);\chiz)\bigr)
 \qquad(q\to\infty),
\end{equation}
where
\begin{equation}\label{eq:Y}
 Y_A(x,q)=\bigl(\log q+(\log x)^2\bigr)(\log_2 q)^A.
\end{equation}
\end{conjecture}

The two terms in \(Y_A(x,q)\) describe different length scales.
For \(x=(\log q)^B\), with \(B>0\) fixed, one has
\(Y_A(x,q)=(\log q)^{1+o(1)}\). For \(x=q^\theta\), with
\(\theta>0\) fixed, the same parameter has size
\((\log q)^2(\log_2q)^A\). The conjecture therefore specifies
how far the Euler product should be truncated as the length of the
sum changes.

No primitivity assumption is made in Conjecture~\ref{conj:GS}.
A fixed primitive character may consequently be induced to larger
moduli, introducing zeros at additional primes while preserving its
values on integers coprime to the new modulus. If the new modulus
contains every prime up to \(z\), then at the cutoff \(z^2\)
the surviving integers, apart from \(1\), are primes.
This observation converts a bound for an ordinary character sum
into cancellation in a prime sum for the fixed inducing character.
Only an upper-bound consequence of the conjecture is needed.

\begin{theorem}\label{thm:implication}
Suppose that there are fixed constants \(A>0\) and \(K>0\) such that,
for every sufficiently large \(q\), every nonprincipal character
\(\chi\bmod q\), and every \(1\le x\le q\),
\begin{equation}\label{eq:upperhyp}
 |S(x,\chi)|
 \le K\Psi(x,Y_A(x,q);\chiz).
\end{equation}
Then the Riemann hypothesis holds for every primitive nonprincipal
Dirichlet \(L\)-function.
\end{theorem}

\begin{corollary}\label{cor:GSrh}
Conjecture~\ref{conj:GS} implies the Riemann hypothesis for every
primitive nonprincipal Dirichlet \(L\)-function.
The same conclusion follows if the error term in
\eqref{eq:GSconj} is replaced by
\(O(\Psi(x,Y_A(x,q);\chiz))\), uniformly in the same range.
\end{corollary}

Indeed, the elementary inequality
\[
 |\Psi(x,y;\chi)|\le\Psi(x,y;\chiz)
\]
shows that either approximation in the corollary implies
\eqref{eq:upperhyp} with a fixed \(K\).

The quantification over imprimitive characters is essential to this
argument. The modulus varies while the inducing primitive character
remains fixed, so the resulting prime-sum estimate controls one
fixed \(L\)-function throughout the half-plane
\(\Rea s>\tfrac{1}{2}\). This gives a different use of character-sum
information from ordinary cancellation estimates of the form
\(S(x,\chi)=o(x)\), whose relationship with the distribution of
zeros is investigated in \cite{GSzeros}.

In the reverse direction, Montgomery and Vaughan's work
\cite{MV} introduced a conditional approximation by sums over
smooth integers. Their argument connects the zero-free region
provided by GRH with a short Euler-product representation.
Granville and Soundararajan refined this method to obtain an error
measured against the smooth-number count, rather than the full
length of the original sum. More precisely, Granville and Soundararajan
\cite[Theorem~2]{GS} showed, under the Riemann hypothesis for
\(L(s,\chi)\), that
\begin{equation}\label{eq:GSthm}
 S(x,\chi)=\Psi(x,y;\chi)
   +O\left(\frac{\Psi(x,y)}{(\log_2 q)^2}\right)
\end{equation}
for \(1\le x\le q\) and
\[
 y\ge(\log q)^2(\log x)^2(\log_2 q)^{12}.
\]
They also obtained a uniform upper bound in terms of
\(\Psi(x,(\log q)^2(\log_2 q)^{20})\). These estimates have subsequent
applications to extreme values of derivatives of Dirichlet \(L\)-functions;
see, for example, Yang \cite{Yang}. In such applications,
partial summation transfers information about the length and size
of character sums to weighted Dirichlet polynomials. The dependence
on the smoothness parameter can therefore affect the constants
obtained for extreme values, making the approximation itself a
useful object of study.

The leading constant in an upper bound for \(L(1,\chi)\)
already illustrates the relevance of the smoothness scale.
Write \(\gamma_{\mathrm E}\) for the Euler--Mascheroni constant.
For nonprincipal characters modulo \(q\), Littlewood's bound
under GRH \cite{Littlewood} is
\[
 |L(1,\chi)|\le
 \bigl(2\e^{\gamma_{\mathrm E}}+o(1)\bigr)\log_2q.
\]
The full Granville--Soundararajan conjecture yields instead
\[
 |L(1,\chi)|\le
 \bigl(\e^{\gamma_{\mathrm E}}+o(1)\bigr)\log_2q.
\]
Both statements are uniform over nonprincipal characters and
apply to general moduli, including imprimitive characters;
see the Introduction of the published version of \cite{Yang}.
These are existing consequences, recalled here to explain the
importance of the cutoff. The fixed-power result below does
not establish the smaller leading constant under GRH.

For sums whose length is at least a fixed power of the modulus, the
approximation can be obtained with the parameter in
\eqref{eq:Y}. Here \(\log x\asymp\log q\), and the proposed
threshold has size \((\log q)^2\) up to powers of
\(\log_2q\). In this range we estimate the ordinary sum and its
smooth restriction separately, showing that both are small compared
with the principal smooth-number count. Their difference is then
small on the same scale. The following form retains a quantitative
error.

\begin{theorem}\label{thm:reverse}
Fix \(0<a\le1\), \(C>6\), and
\[
 0<\gamma<\frac{C-6}{4}.
\]
Let \(\chi\bmod q\) be nonprincipal, and suppose that the Riemann
hypothesis holds for \(L(s,\chi)\). Uniformly for
\[
 q^a\le x\le q,\qquad y=Y_C(x,q),
\]
we have
\begin{equation}\label{eq:strong}
 \begin{split}
 |S(x,\chi)|+|\Psi(x,y;\chi)|
 &\ll_{a,C,\gamma}\Psi(x,y;\chiz)\\
 &\quad{}\times
 \exp\left\{-\gamma
       \frac{(\log x)\log_3 q}{\log_2 q}\right\}.
 \end{split}
\end{equation}
Consequently,
\begin{equation}\label{eq:reverseapprox}
 S(x,\chi)=\Psi(x,y;\chi)
          +o_{a,C}\bigl(\Psi(x,y;\chiz)\bigr).
\end{equation}
The estimates apply to arbitrary moduli and include imprimitive
characters.
\end{theorem}

The normalization in \eqref{eq:reverseapprox} is the nonnegative
count \(\Psi(x,y;\chiz)\). No relative-error assertion with
denominator \(|\Psi(x,y;\chi)|\) is intended; the smooth character
sum may itself exhibit substantial cancellation. The constants
may depend on \(a\), and no uniformity as \(a\to0\) is asserted.

For example, \(C=8\) is admissible. In the range \(q^a\le x\le q\),
the parameter in \eqref{eq:GSthm} has order
\((\log q)^4(\log_2 q)^{12}\), whereas the choice \(Y_8(x,q)\)
has order \((\log q)^2(\log_2 q)^8\). Thus
Theorem~\ref{thm:reverse} reduces the required power of \(\log q\)
for this part of the approximation.

The mechanism in Theorem~\ref{thm:reverse} is that both character sums
are small compared with \(\Psi(x,y;\chiz)\). The argument does not
require an asymptotic formula for the latter: products of distinct primes
in \((\tfrac{y}{2},y]\) provide a sufficient lower bound. The number of
primes excluded by divisibility by \(q\) is negligible relative
to the available primes in this interval. This yields the required
uniformity over arbitrary moduli without replacing
\(\Psi(x,y;\chiz)\) by an unrestricted smooth-number count.

On the analytic side, bounds for \(\log L(s,\chi)\) and for the
logarithm of its finite Euler product allow the Perron contour to
be moved to
\[
 \Rea s=\frac12+b\frac{\log_3q}{\log_2q}.
\]
The analytic error is small enough when \(b>\tfrac{3}{2}\), while the
combinatorial lower bound provides a coefficient \(\tfrac{C}{4}\) on
the scale \((\log x)\log_3q/\log_2q\). The choice
\(\tfrac{3}{2}<b<\tfrac{C}{4}\) explains the condition \(C>6\) in the theorem.
The quantitative saving follows by retaining a positive gap between
these two coefficients. The condition \(C>6\) is a sufficient
condition for this argument; no optimality is claimed.

\subsection*{Zeta sums}

We now state the zeta-sum counterparts of the character-sum results.

For real \(t\), define
\[
 S(x,t)=\sum_{n\le x}n^{-\ii t},\qquad
 \Psi(x,y;t)=
       \sum_{\substack{n\le x\\\PP(n)\le y}}n^{-\ii t}.
\]
Yang \cite{Yang} formulated the following analogue of
Conjecture~\ref{conj:GS}.\footnote{Numbered results of Yang cited in
this paper refer to \texttt{arXiv:2204.13826v4}.}

\begin{conjecture}[Yang, Conjecture~1]\label{conj:zeta}
There is a constant \(A>0\) such that, uniformly for
\[
 1\le x\le T,\qquad 2T\le t\le5T,
\]
we have
\[
 S(x,t)=\Psi(x,Y_A(x,T);t)
            +o\bigl(\Psi(x,Y_A(x,T))\bigr)
 \qquad(T\to\infty).
\]
\end{conjecture}

Under RH, Yang \cite[Theorem~5]{Yang} established a smooth
approximation in the range
\[
 2\le x\le T,\qquad T+y+3\le t\le T^{1000},
\]
provided that
\[
 y\ge(\log T)^2(\log x)^2(\log_2T)^{12}.
\]
For the polylogarithmic cutoffs considered below, the height
restrictions hold throughout \(2T\le t\le5T\) once \(T\) is
sufficiently large. In the range \(x\ge T^a\), we have the
following counterpart of Theorem~\ref{thm:reverse}.

\begin{theorem}\label{thm:zeta}
Assume the Riemann hypothesis for \(\zeta(s)\).
Fix \(0<a\le1\), \(C>6\), and \(0<\gamma<\tfrac{C-6}{4}\).
Uniformly for
\[
 T^a\le x\le T,\qquad 2T\le t\le5T,\qquad y=Y_C(x,T),
\]
we have
\begin{equation}\label{eq:zetastrong}
 \begin{split}
 |S(x,t)|+|\Psi(x,y;t)|
 &\ll_{a,C,\gamma}\Psi(x,y)\\
 &\quad{}\times
 \exp\left\{-\gamma
       \frac{(\log x)\log_3T}{\log_2T}\right\}.
 \end{split}
\end{equation}
In particular,
\[
 S(x,t)=\Psi(x,y;t)+o_{a,C}\bigl(\Psi(x,y)\bigr).
\]
\end{theorem}

The full zeta-sum conjecture also has a consequence for values on the
critical line. The following result requires only its upper-bound
consequence.

\begin{proposition}\label{prop:lindelof}
Conjecture~\ref{conj:zeta} implies the Lindel\"of hypothesis:
for every \(\varepsilon>0\),
\[
 \zeta(\tfrac{1}{2}+\ii t)\ll_\varepsilon (1+|t|)^\varepsilon.
\]
In fact, the uniform upper bound
\[
 |S(x,t)|\ll\Psi(x,Y_A(x,T)),
 \qquad 1\le x\le T,\quad 2T\le t\le5T,
\]
for some fixed \(A>0\), suffices.
\end{proposition}

The proof combines Rankin's inequality with the approximate functional
equation. Its conclusion is the Lindel\"of hypothesis; no implication
from the zeta-sum conjecture to RH is asserted here.

\begin{remark}\label{rem:scope}
Throughout this note, a Riemann-hypothesis assertion about a
nonprincipal character refers to its associated Dirichlet \(L\)-function.
For an imprimitive character this is equivalent to the corresponding
assertion for its inducing primitive character, since the additional
Euler factors have zeros only on \(\Rea s=0\).
Theorem~\ref{thm:implication} and Corollary~\ref{cor:GSrh} make no
assertion about the Riemann hypothesis for \(\zeta(s)\): the hypothesis
excludes principal characters.
\end{remark}

The paper is organized as follows. Section~2 proves
Theorem~\ref{thm:implication} by induction to moduli containing all
small primes. Section~3 develops the logarithmic and Perron
estimates, and Section~4 gives the lower bound for smooth integers
coprime to the modulus. These estimates are combined in Section~5
to prove Theorem~\ref{thm:reverse}. Section~\ref{sec:zeta}
contains the proofs of Theorem~\ref{thm:zeta} and
Proposition~\ref{prop:lindelof}.

\section{Induction to moduli containing all small primes}

\begin{proof}[Proof of Theorem~\ref{thm:implication}]
Fix a primitive nonprincipal character \(\xi\bmod r\).
For a sufficiently large real number \(z>r\), put
\begin{equation}\label{eq:primorial}
 q_z=r\prod_{p\le z}p,\qquad X=z^2,
\end{equation}
and let \(\chi_z\bmod q_z\) be induced by \(\xi\). Thus
\[
 \chi_z(n)=
 \begin{cases}
  \xi(n),&(n,q_z)=1,\\
  0,&(n,q_z)>1.
 \end{cases}
\]
The reduction map
\((\mathbb Z/q_z\mathbb Z)^\times\to
(\mathbb Z/r\mathbb Z)^\times\) is surjective, so \(\chi_z\)
is nonprincipal. The prime number theorem gives
\[
 \log q_z=\log r+\sum_{p\le z}\log p\sim z.
\]
In particular, \(X\le q_z\) for large \(z\). Write
\[
 y_z=Y_A(X,q_z).
\]
Since \(\log X=2\log z\), we have
\begin{equation}\label{eq:yz}
 y_z\sim z(\log z)^A,\qquad z<y_z<X.
\end{equation}

Every prime divisor of an integer coprime to \(q_z\) exceeds \(z\).
If such an integer is composite, it is larger than \(z^2\).
It follows that
\begin{equation}\label{eq:survive}
 S(X,\chi_z)=1+\sum_{z<p\le X}\xi(p).
\end{equation}
Since all prime divisors of \(r\) are at most \(z\), we also have
\begin{equation}\label{eq:principalcount}
 \Psi(X,y_z;\chi_{0,z})=1+\pi(y_z)-\pi(z),
\end{equation}
where \(\chi_{0,z}\) denotes the principal character modulo \(q_z\).
Using \eqref{eq:upperhyp}, \eqref{eq:survive}, and
\eqref{eq:principalcount}, we obtain
\[
 \left|\sum_{z<p\le X}\xi(p)\right|
 \ll 1+\pi(y_z)
 \ll \frac{y_z}{\log y_z}.
\]
Adding the primes at most \(z\), and using \eqref{eq:yz}, gives
\[
 \left|\sum_{p\le X}\xi(p)\right|
 \ll_{\xi,A,K} z(\log z)^{A-1}.
\]
Here \(\pi(z)\ll z/\log z\) is absorbed because \(A>0\).
As \(z=\sqrt X\) is an arbitrary sufficiently large real number,
\begin{equation}\label{eq:primebound}
 P_\xi(X):=\sum_{p\le X}\xi(p)
 \ll_{\xi,A,K}X^{\frac{1}{2}}(\log X)^{A-1}
\end{equation}
holds for all sufficiently large \(X\).

Partial summation now yields
\begin{align}
 \vartheta(X,\xi)
 &:=\sum_{p\le X}\xi(p)\log p \notag\\
 &=P_\xi(X)\log X-\int_2^X\frac{P_\xi(u)}{u}\dd u
 \ll_{\xi,A,K}X^{\frac{1}{2}}(\log X)^A.
 \label{eq:theta}
\end{align}
For the prime powers, Chebyshev's estimate gives
\[
 \sum_{\substack{p^j\le X\\j\ge2}}\log p
 \ll X^{\frac{1}{2}}+X^{\frac{1}{3}}\log X
 \ll X^{\frac{1}{2}}.
\]
Consequently,
\begin{equation}\label{eq:psi}
 \psi(X,\xi):=\sum_{n\le X}\Lambda(n)\xi(n)
 \ll_{\xi,A,K}X^{\frac{1}{2}}(\log X)^A.
\end{equation}

The integral
\[
 H_\xi(s)=s\int_1^\infty\psi(u,\xi)u^{-s-1}\dd u
\]
converges absolutely and locally uniformly on \(\Rea s>\tfrac{1}{2}\),
by \eqref{eq:psi}. Hence it defines a holomorphic function there.
For \(\Rea s>1\), summation by parts gives
\[
 H_\xi(s)=\sum_{n=1}^\infty\frac{\Lambda(n)\xi(n)}{n^s}
         =-\frac{L'}{L}(s,\xi).
\]
Since \(L(s,\xi)\) is entire, the function
\(L'(s,\xi)+H_\xi(s)L(s,\xi)\) is holomorphic on
\(\Rea s>\tfrac{1}{2}\). It vanishes on \(\Rea s>1\), and the
identity theorem shows that it vanishes throughout that half-plane.
If \(L(s,\xi)\) had a zero of order \(m\ge1\) there, then
\(L'(s,\xi)\) would have order \(m-1\), whereas
\(H_\xi(s)L(s,\xi)\) would have order at least \(m\), a contradiction.
Thus \(L(s,\xi)\) has no zero in \(\Rea s>\tfrac{1}{2}\).
For a nontrivial zero \(\rho\), the functional equation and complex
conjugation show that \(1-\overline{\rho}\) is also a zero of
\(L(s,\xi)\). This excludes \(\Rea\rho<\tfrac{1}{2}\) as well.
\end{proof}

\begin{remark}
Only the family \eqref{eq:primorial} is used for each fixed \(\xi\).
Neither the conjectured little-\(o\) error nor uniformity in the
primitive conductor is needed for this implication once
\eqref{eq:upperhyp} is available along that family.
A version restricted to primitive characters or to prime moduli
would not supply these induced characters.
\end{remark}

\section{Bounds for a Dirichlet \texorpdfstring{\(L\)}{L}-function and its finite Euler product}

We prepare the estimates needed for Theorem~\ref{thm:reverse}.
Write
\[
 Q=\log q,\qquad \ell=\log Q.
\]
For a nonprincipal character \(\chi\bmod q\), define
\begin{equation}\label{eq:F}
 F_y(s,\chi)=\prod_{p\le y}
        \left(1-\frac{\chi(p)}{p^s}\right)^{-1}.
\end{equation}
This is holomorphic and nonzero on \(\Rea s>0\), where its
absolutely convergent Dirichlet series is
\[
 F_y(s,\chi)=\sum_{\PP(n)\le y}\frac{\chi(n)}{n^s}.
\]
Under the Riemann hypothesis for \(L(s,\chi)\), we use the branch of
\(\log L(s,\chi)\) on \(\Rea s>\tfrac{1}{2}\) obtained by continuation of its
Euler-product logarithm from \(\Rea s>1\).
The logarithm of \(F_y\) is defined by its absolutely convergent
prime-power expansion on \(\Rea s>0\).

\begin{lemma}\label{lem:log}
Fix \(b>0\), \(C>0\), and positive constants \(d_1,d_2\).
Suppose that
\[
 d_1Q^2\ell^C\le y\le d_2Q^2\ell^C,
 \qquad
 \sigma_b=\frac12+b\frac{\log\ell}{\ell}.
\]
If the Riemann hypothesis holds for \(L(s,\chi)\), then, for all
sufficiently large \(q\),
\begin{equation}\label{eq:logbound}
 |\log L(s,\chi)|+|\log F_y(s,\chi)|
 \ll_{b,C,d_1,d_2}Q\ell^{2-2b}
\end{equation}
uniformly for \(\sigma_b\le\Rea s\le2\) and \(|\Ima s|\le q\).
\end{lemma}

\begin{proof}
For \(\Rea s\ge1+1/Q\), the Euler-product logarithms give
\[
 |\log L(s,\chi)|+|\log F_y(s,\chi)|
 \le2\log\zeta(1+1/Q)\ll\ell,
\]
which is absorbed by the right-hand side of \eqref{eq:logbound}.
We may therefore suppose that \(\Rea s\le1+1/Q\).
We use the logarithmic approximation in
\cite[\S8, Lemma~8.2]{GS}; the lemma number here and below is that
of \texttt{arXiv:math/9903196v1}.
For \(v\) of order \(Q^2\) or \(Q^2\ell^C\), we have
\(v=o(q)\). Thus the auxiliary heights satisfy
\(|\Ima s|+v+3\le q+v+3<3q\) for large \(q\), as required
by the cited estimates. Write \(\sigma=\Rea s\).
The error term in that lemma is
\[
 O\left(\frac{Q}{(\sigma_1-\sigma_0)^2}
             v^{\sigma_1-\sigma}\right),\qquad
 \sigma_1=\min\left\{\frac{\sigma+\sigma_0}{2},
                      \sigma_0+\frac1{\log v}\right\}.
\]
To keep the zero-free boundary strictly to the right of the
critical line, choose
\[
 \sigma_0=\frac12+\frac1{\log v}.
\]
For sufficiently large \(q\), we have
\(\sigma-\tfrac12\ge b(\log\ell)/\ell>3/\log v\), so
\(\sigma_1=\tfrac12+2/\log v\).
The zero-free hypothesis follows from the Riemann hypothesis.
Since \(\sigma_1-\sigma_0=1/\log v\) and
\(v^{2/\log v}=\e^2\), we obtain
\begin{equation}\label{eq:logapprox}
 \log L(s,\chi)
 =\sum_{2\le n\le v}
      \frac{\Lambda(n)\chi(n)}{n^s\log n}
  +O\left(Q(\log v)^2v^{\frac{1}{2}-\Rea s}\right).
\end{equation}
One may initially take \(v\) to be a half-integer; replacing a real
cutoff by a neighbouring half-integer affects the displayed sum by at
most one term and does not alter the estimates below.

Apply \eqref{eq:logapprox} with \(v\) a neighbouring half-integer to
\(Q^2\). Since the absolute value of each coefficient
\(\Lambda(n)/\log n\) is at most \(1\),
\[
 \left|\sum_{2\le n\le v}
      \frac{\Lambda(n)\chi(n)}{n^s\log n}\right|
 \le\sum_{n\le v}n^{-\sigma_b}
 \ll Q\ell^{-2b}.
\]
Also,
\[
 Q(\log v)^2v^{\frac{1}{2}-\sigma_b}
 \ll Q\ell^{2-2b}.
\]
This proves the asserted bound for \(\log L\).

Next, comparison of prime-power expansions gives
\begin{equation}\label{eq:Fcomparison}
 \log F_y(s,\chi)
 =\sum_{2\le n\le y}
      \frac{\Lambda(n)\chi(n)}{n^s\log n}
 +O\left(\sum_p\sum_{j\ge2}p^{-j\Rea s}\right).
\end{equation}
If \(\eta=\sigma_b-\tfrac12\), the double sum is
\[
 \ll\sum_{n\ge2}n^{-1-2\eta}\ll\eta^{-1}.
\]
By \eqref{eq:logapprox} at \(v=y\), the difference between the
main sum in \eqref{eq:Fcomparison} and \(\log L(s,\chi)\) is
\[
 \ll Q(\log y)^2y^{\frac{1}{2}-\sigma_b}
 \ll Q\ell^{2-2b},
\]
because \(\log y=2\ell+C\log\ell+O(1)\) and \(C>0\).
The term \(\eta^{-1}\ll_b\ell/\log\ell\) is absorbed by
\(Q\ell^{2-2b}\). This proves \eqref{eq:logbound}.
\end{proof}

\begin{lemma}\label{lem:perron}
Under the assumptions of Lemma~\ref{lem:log}, fix \(0<a\le1\).
Uniformly for \(q^a\le x\le q\), we have
\begin{equation}\label{eq:perronbound}
 \begin{split}
 \log\bigl(1+|S(x,\chi)|+|\Psi(x,y;\chi)|\bigr)
 \le {}&\frac{\log x}{2}
       +b\frac{(\log x)\log\ell}{\ell}\\
       &+O_{a,b,C,d_1,d_2}
         \bigl(Q\ell^{2-2b}+\ell\bigr).
 \end{split}
\end{equation}
\end{lemma}

\begin{proof}
Put \(z=\lfloor x\rfloor+\tfrac12\), keeping \(y\) fixed.
The sums with cutoff \(x\) equal those with cutoff \(z\).
Consider either
\[
 \mathcal A(s)=L(s,\chi)
 \quad\hbox{or}\quad
 \mathcal A(s)=F_y(s,\chi).
\]
In both cases the coefficients \(a_n\) of the Dirichlet series satisfy
\(|a_n|\le1\). With \(c=1+1/\log z\) and \(T=q\), truncated Perron
inversion gives
\begin{equation}\label{eq:truncperron}
 \sum_{n\le z}a_n
 =\frac1{2\pi\ii}\int_{c-\ii q}^{c+\ii q}
       \mathcal A(s)\frac{z^s}{s}\dd s
 +O(\log^2q).
\end{equation}
For completeness, the usual error is bounded by
\[
 \sum_{n\ge1}\left(\frac zn\right)^c
 \min\left\{1,\frac1{q|\log(z/n)|}\right\}.
\]
For \(\tfrac{z}{2}<n<2z\), use
\(|\log(z/n)|\asymp|z-n|/z\) and \(|z-n|\ge\tfrac12\).
The remaining ranges contribute
\(O(z^c\zeta(c)/q)\). These bounds justify the error in
\eqref{eq:truncperron}, since \(z\le q+\tfrac12\).

Let \(M\) be a sufficiently large constant multiple of
\(Q\ell^{2-2b}\). Lemma~\ref{lem:log} gives
\[
 |\mathcal A(s)|\le\e^M
 \qquad
 (\sigma_b\le\Rea s\le c,\ |\Ima s|\le q).
\]
Neither of the two functions has a pole in this rectangle, so the contour
may be shifted to \(\Rea s=\sigma_b\). The two horizontal
segments contribute \(O(\e^M z^c/q)=O(\e^M)\).
The new vertical segment contributes
\[
 \ll \e^M z^{\sigma_b}\int_{-q}^q
                     \frac{\dd t}{|\sigma_b+\ii t|}
 \ll \e^M z^{\sigma_b}\log q.
\]
Using \(\log z=\log x+O(1/x)\) and taking logarithms proves
\eqref{eq:perronbound}.
\end{proof}

\section{A lower bound for smooth integers coprime to the modulus}

\begin{lemma}\label{lem:lower}
Fix \(0<a\le1\), \(C>0\), and positive constants \(d_1,d_2\).
Uniformly for
\[
 q^a\le x\le q,\qquad
 d_1Q^2\ell^C\le y\le d_2Q^2\ell^C,
\]
we have
\begin{equation}\label{eq:lower}
 \log\Psi(x,y;\chiz)
 \ge \frac{\log x}{2}
  +\left(\frac C4+o_{a,C,d_1,d_2}(1)\right)
       \frac{(\log x)\log\ell}{\ell}.
\end{equation}
\end{lemma}

\begin{proof}
Let
\[
 \mathcal P=\{p:\tfrac{y}{2}<p\le y,\ p\nmid q\},\qquad N=|\mathcal P|.
\]
The number of excluded primes satisfies
\[
 0\le \pi(y)-\pi(\tfrac{y}{2})-N\le\omega(q)\le Q/\log2.
\]
Consequently, the prime number theorem gives
\begin{equation}\label{eq:N}
 N\sim\frac{y}{2\log y},
\end{equation}
uniformly in \(q\), since
\(Q=o(y/\log y)\).
Write \(X=\log x\), \(L=\log y\), and
\[
 m=\left\lfloor\frac XL\right\rfloor.
\]
Products of \(m\) distinct members of \(\mathcal P\) are at most
\(y^m\le x\), are \(y\)-smooth, and are coprime to \(q\).
Unique factorization therefore gives
\begin{equation}\label{eq:binomial}
 \Psi(x,y;\chiz)\ge \binom{N}{m}
                  \ge\left(\frac Nm\right)^m.
\end{equation}
Here \(m\to\infty\) and \(N/m\to\infty\) uniformly in the stated
range.

We have
\[
 L=2\ell+C\log\ell+O_{C,d_1,d_2}(1),\qquad
 \log X=\ell+O_a(1).
\]
It follows from \eqref{eq:N} that
\[
 \log(N/m)=L-\log X-\log2+o(1)
          =\ell+C\log\ell+O_{a,C,d_1,d_2}(1).
\]
Taking logarithms in \eqref{eq:binomial} now gives
\begin{align*}
 \log\Psi(x,y;\chiz)
 &\ge
 \frac{X\bigl(\ell+C\log\ell+O(1)\bigr)}
      {2\ell+C\log\ell+O(1)}
 -O(\ell)\\
 &\ge \frac X2+\frac C4\frac{X\log\ell}{\ell}
 -O_{a,C,d_1,d_2}\left(
     \frac{X(\log\ell)^2}{\ell^2}
       +\frac X\ell+\ell\right).
\end{align*}
The error is \(o(X\log\ell/\ell)\), which proves
\eqref{eq:lower}.
\end{proof}

\section{Proof of the conditional approximation}

\begin{proof}[Proof of Theorem~\ref{thm:reverse}]
Fix \(a,C,\gamma\) as in the theorem, and choose \(b\) with
\begin{equation}\label{eq:bchoice}
 \frac32<b<\frac C4-\gamma.
\end{equation}
Such a choice is possible because \(\gamma<\tfrac{C-6}{4}\).
Put \(Q=\log q\), \(\ell=\log Q\), \(X=\log x\), and
\[
 y=(Q+X^2)\ell^C.
\]
Since \(aQ\le X\le Q\), we have, for large \(q\),
\[
 a^2Q^2\ell^C\le y\le2Q^2\ell^C.
\]
Thus Lemmas~\ref{lem:log}--\ref{lem:lower} all apply.

Let
\[
 W=\frac{X\log\ell}{\ell}.
\]
Lemma~\ref{lem:perron} gives
\[
 \log\bigl(1+|S(x,\chi)|+|\Psi(x,y;\chi)|\bigr)
 \le \frac X2+bW+O(Q\ell^{2-2b}+\ell).
\]
By \(X\ge aQ\) and \(b>\tfrac{3}{2}\),
\[
 \frac{Q\ell^{2-2b}+\ell}{W}
 \ll_a
 \frac{\ell^{3-2b}}{\log\ell}
 +\frac{\ell^2}{Q\log\ell}
 =o(1).
\]
On the other hand, Lemma~\ref{lem:lower} shows that
\[
 \log\Psi(x,y;\chiz)\ge\frac X2+
                         \left(\frac C4+o(1)\right)W.
\]
Subtracting these inequalities yields
\[
 \log\frac{1+|S(x,\chi)|+|\Psi(x,y;\chi)|}
                {\Psi(x,y;\chiz)}
 \le-\left(\frac C4-b+o(1)\right)W.
\]
The strict inequality in \eqref{eq:bchoice} proves
\eqref{eq:strong} for sufficiently large \(q\).
Finally,
\[
 |S(x,\chi)-\Psi(x,y;\chi)|
 \le |S(x,\chi)|+|\Psi(x,y;\chi)|
\]
gives \eqref{eq:reverseapprox}.
\end{proof}

\begin{corollary}\label{cor:eight}
Assume the Riemann hypothesis for \(L(s,\chi)\), where
\(\chi\bmod q\) is nonprincipal. Fix \(0<a\le1\), and put
\[
 y=\bigl(\log q+(\log x)^2\bigr)(\log_2 q)^8.
\]
Uniformly for \(q^a\le x\le q\),
\[
 S(x,\chi)=\Psi(x,y;\chi)
 +O_a\left(
   \Psi(x,y;\chiz)
   \exp\left\{-\frac14
        \frac{(\log x)\log_3 q}{\log_2 q}\right\}\right).
\]
\end{corollary}

\begin{proof}
Take \(C=8\) and \(\gamma=\tfrac{1}{4}\) in
Theorem~\ref{thm:reverse}.
\end{proof}

\section{Proofs for zeta sums}\label{sec:zeta}

We prove Theorem~\ref{thm:zeta} and Proposition~\ref{prop:lindelof},
using the notation and hypotheses stated in the Introduction.

\begin{proof}[Proof of Theorem~\ref{thm:zeta}]
Put
\[
 Q=\log T,\quad \ell=\log Q,\quad X=\log x,\quad
 \sigma_b=\frac12+b\frac{\log\ell}{\ell},
\]
where \(\tfrac{3}{2}<b<\tfrac{C}{4}-\gamma\).
We have \(y=(Q+X^2)\ell^C\asymp_a Q^2\ell^C\).
Write
\[
 F_y(w)=\prod_{p\le y}(1-p^{-w})^{-1}.
\]
Under RH, \(\zeta(w)\) is holomorphic and nonzero on the
simply connected domain \(\Rea w>\tfrac{1}{2},\ \Ima w>0\).
We choose \(\log\zeta(w)\) by continuation of the
Euler-product logarithm from \(\Rea w>1\).
The logarithm of \(F_y(w)\) is defined by its absolutely
convergent prime-power expansion on \(\Rea w>0\).
The logarithmic approximation
\cite[Lemma~3]{Yang} gives, for \(v\asymp Q^2\) or
\(v\asymp Q^2\ell^C\),
\begin{equation}\label{eq:zetalogapprox}
 \log\zeta(\sigma+\ii\tau)
 =\sum_{2\le n\le v}\frac{\Lambda(n)}
                {n^{\sigma+\ii\tau}\log n}
 +O\left(Q(\log v)^2v^{\frac{1}{2}-\sigma}\right)
\end{equation}
uniformly for \(\sigma_b\le\sigma\le1\) and
\(T\le\tau\le6T\).
Indeed, its hypotheses \(\tau\ge v+3\) hold for large \(T\);
the minimum defining the auxiliary abscissa in that lemma is
\(\tfrac12+1/\log v\) in this range.

Taking \(v\asymp Q^2\) bounds the finite sum by
\(O(Q\ell^{-2b})\) and the error by
\(O(Q\ell^{2-2b})\). Taking \(v=y\), and comparing the
prime-power series as in \eqref{eq:Fcomparison}, then yields
\begin{equation}\label{eq:zetalogbounds}
 |\log\zeta(\sigma+\ii\tau)|+
 |\log F_y(\sigma+\ii\tau)|
 \ll_{a,b,C}Q\ell^{2-2b}.
\end{equation}
Here the estimate is uniform for \(\sigma_b\le\sigma\le1\)
and \(T\le\tau\le6T\).
The prime-power correction is
\(O((\sigma_b-\tfrac12)^{-1})\) and is absorbed by this bound.
For \(1\le\sigma\le2\), the elementary estimates
\[
 |\zeta(\sigma+\ii\tau)|\ll\log(2+\tau),\qquad
 |F_y(\sigma+\ii\tau)|\le F_y(1)\ll\log y
\]
supply the corresponding modulus bounds in the remaining part of
the Perron rectangle. The zeta bound follows, for example, by
Euler--Maclaurin truncation at a point of order \(\tau\).
Combining these estimates, there is a constant \(K=K(a,b,C)>0\)
such that
\[
 |\zeta(w)|+|F_y(w)|
 \le \exp\{KQ\ell^{2-2b}\}
\]
throughout
\(\sigma_b\le\Rea w\le2,\ T\le\Ima w\le6T\),
for sufficiently large \(T\). These modulus bounds justify
the contour shift below.

For fixed \(t\in[2T,5T]\), use truncated Perron inversion at
height \(T\) for the Dirichlet series
\[
 \mathcal A(s)=\zeta(s+\ii t)
 \quad\hbox{and}\quad
 \mathcal A(s)=F_y(s+\ii t).
\]
Their coefficients have absolute value at most \(1\).
As in Lemma~\ref{lem:perron}, replace the cutoff \(x\) by
\(z=\lfloor x\rfloor+\tfrac12\), keeping \(y\) fixed.
The truncation error is \(O(\log^2T)\).
On the rectangle with \(|\Ima s|\le T\), the imaginary part
of \(s+\ii t\) belongs to \([T,6T]\). In particular, the pole
of \(\zeta(s+\ii t)\), at \(s=1-\ii t\), is outside the
rectangle. The contour can therefore be shifted to
\(\Rea s=\sigma_b\) without a residue.

The same vertical and horizontal estimates used in
Lemma~\ref{lem:perron}, now with \eqref{eq:zetalogbounds},
give
\begin{equation}\label{eq:zetaperron}
 \log\bigl(1+|S(x,t)|+|\Psi(x,y;t)|\bigr)
 \le \frac X2+b\frac{X\log\ell}{\ell}
       +O_{a,b,C}(Q\ell^{2-2b}+\ell).
\end{equation}
The proof of Lemma~\ref{lem:lower}, with all primes in
\((\tfrac{y}{2},y]\) available, gives
\[
 \log\Psi(x,y)\ge \frac X2+
       \left(\frac C4+o_{a,C}(1)\right)
       \frac{X\log\ell}{\ell}.
\]
Since \(X\ge aQ\) and \(b>\tfrac{3}{2}\), the error in
\eqref{eq:zetaperron} is \(o(X\log\ell/\ell)\).
The comparison made in the proof of
Theorem~\ref{thm:reverse} proves \eqref{eq:zetastrong}.
\end{proof}

\begin{proof}[Proof of Proposition~\ref{prop:lindelof}]
Take \(T\) large and \(t\in[2T,5T]\). For \(1\le u\le T\),
\[
 Y_A(u,T)\le Y:=2(\log T)^2(\log_2T)^A.
\]
Fix \(0<\delta<\tfrac{1}{2}\), and put \(\alpha=\tfrac{1}{2}+\delta\).
Rankin's inequality gives
\[
 \Psi(u,Y)\le u^\alpha
       \prod_{p\le Y}(1-p^{-\alpha})^{-1}.
\]
The logarithm of the product is \(O_\delta(Y^{1-\alpha})\), and
\[
 Y^{1-\alpha}
 \ll_{\delta,A}
 (\log T)^{1-2\delta}(\log_2T)^{A(\frac{1}{2}-\delta)}
 =o(\log T).
\]
The assumed upper bound therefore yields
\begin{equation}\label{eq:zetaLHsum}
 |S(u,t)|\ll_{\delta,A}u^{\frac{1}{2}+\delta}T^{o(1)}
 \qquad(1\le u\le T).
\end{equation}
The factor \(T^{o(1)}\) here is uniform in \(u\) and \(t\)
in the stated ranges, with \(\delta,A\) fixed.
Let \(N=\lfloor\sqrt{t/(2\pi)}\rfloor\); then \(N\le T\)
for large \(T\). Partial summation gives
\[
 \begin{split}
 \sum_{n\le N}n^{-\frac{1}{2}-\ii t}
 &=N^{-\frac{1}{2}}S(N,t)
       +\frac12\int_1^N S(u,t)u^{-\frac{3}{2}}\dd u\\
 &\ll_{\delta,A}T^{o(1)}N^\delta
 \ll_{\delta,A}T^{\frac{\delta}{2}+o(1)}.
 \end{split}
\]
The sum with exponent \(-\tfrac{1}{2}+\ii t\) is its complex conjugate.
The approximate functional equation on the critical line
\cite[\S4.17, formula~(4.17.1)]{Titchmarsh} reads
\[
 \zeta(\tfrac{1}{2}+\ii t)
 =\sum_{n\le N}n^{-\frac{1}{2}-\ii t}
   +\mathcal X(t)\sum_{n\le N}n^{-\frac{1}{2}+\ii t}
   +O(t^{-\frac{1}{4}}),
\]
where, with \(\Gamma\) denoting the gamma function,
\[
 \mathcal X(t)=\pi^{\ii t}
       \frac{\Gamma(\tfrac{1}{4}-\tfrac{\ii t}{2})}{\Gamma(\tfrac{1}{4}+\tfrac{\ii t}{2})},
 \qquad |\mathcal X(t)|=1.
\]
Thus \(\zeta(\tfrac{1}{2}+\ii t)\ll_{\delta,A}
T^{\frac{\delta}{2}+o(1)}\).
Given \(\varepsilon>0\), choose
\(0<\delta<\min\{\tfrac{1}{2},\varepsilon\}\).
For sufficiently large \(T\), we then have
\(\tfrac{\delta}{2}+o(1)<\varepsilon\).
For any sufficiently large \(t>0\), take \(T=\tfrac{t}{3}\).
The result follows for positive \(t\), and hence for negative
\(t\) by conjugation. Bounded \(t\) is absorbed into the
implied constant.
\end{proof}

\begin{remark}
The induced-character argument of Section~2 uses vanishing at all
integers divisible by small primes. The factors \(n^{-\ii t}\)
never vanish, so that argument does not transfer to
Conjecture~\ref{conj:zeta}. Proposition~\ref{prop:lindelof}
establishes a Lindel\"of consequence; it makes no claim that the
zeta-sum conjecture implies RH.
\end{remark}
\section*{Acknowledgement}
The authors thank OpenAI's ChatGPT for assistance with the presentation of the manuscript. Zikang Dong is supported by the National Natural Science Foundation of China (Grant No. 1240011770). Guang-Liang Zhou is supported by the National Natural
Science Foundation of China (Grant No. 12401009)

\end{document}